# La variable complexe dans la première décennie de l'oeuvre grothendieckienne : collection d'exemples, réseau d'idées, source de visions

Fernando Zalamea[(*)]

*Abstract*. L'oeuvre de Grothendieck se mesure en termes de souplesse, naturalité, universalité. Ces caractéristiques sont déjà présentes dès ses premiers travaux, en particulier autour de ses contributions à la variable complexe, grande source de souplesse. Nous présentons ici ses exemples et ses idées autour de la variable complexe dans la première décennie (1949-1957), ses références à Riemann, et son influence sur la vision obtenue au ICM 1958 (schémas et topos, en vue de la résolution des conjectures de Weil). Finalement, nous donnons quelques aperçus sur les développements postérieurs (1958-1970 & 1981-1991) de la variable complexe dans l'oeuvre de Grothendieck.

La variable complexe se situe au *coeur* de la pensée mathématique. Le programme de complexification de Riemann[1] a été en effet poussé sur un très grand spectre de domaines mathématiques, et ses effets ont été énormes[2]. Un magnifique diagramme de René Thom pour la *Enciclopedia Einaudi* (1982)[3] montre la "teoria delle funzioni di variabile complessa" au *centre* d'un plan divisé, horizontalement, par la dualité (mathématiques discrètes / mathématiques continues), et, verticalement, par la dualité (génération libre / génération liée). Les fonctions de la variable complexe participent donc, en même temps, du continu et du discret, de la liberté abstraite et des liaisons du concret, ce qui peut expliquer leur souplesse et ductilité. Dans le diagramme de Thom, la variable complexe se projette directement sur les espaces vectoriels topologiques, thème de la Thèse de Grothendieck (1949-53), et se trouve juste à côté de la continuation analytique (donc, des faisceaux) et de la géométrie algébrique, sources du tournant grothendieckien à Kansas (1955). C'est autour de ce *double coeur* formé par la variable complexe —autour des mathématiques et autour de l'émergence du premier Grothendieck— que nous voulons nous arrêter.

---

[(*)] Departamento de Matemáticas, Universidad Nacional de Colombia. https ://unal.academia.edu/FernandoZalamea
[1] Riemann, Bernhard (1851), Thèse Doctorale, "Principes fondamentaux pour une théorie générale des fonctions d'une grandeur variable complexe" (§XX), in *: Oeuvres mathématiques*, Paris : Jacques Gabay, 1990.
[2] Allant, par exemple, de la théorie des nombres (avec la preuve du théorème de Fermat utilisant l'architectonique complexe de Taniyama-Shimura, ou l'hypothèse de Riemann utilisant la continuation analytique) à la physique mathématique (avec le groupe cosmique de Cartier-Connes-Kontsevich, les hypothèses complexes de Penrose, ou la reconstruction par Furey des forces fondamentales avec quaternions et octonions), en passant par toutes sortes d'applications en géométrie algébrique, géométrie différentielle, topologie algébrique, algèbre topologique, analyse fonctionnelle, etc.
[3] Thom, René (1982), "L'aporia fondatrice della matematica", in : *Enciclopedia Einaudi*, Torino : Einaudi, vol. 15, p. 1144.

En quatre sections, notre article étudie : (1) quelques traits généraux –souplesse, naturalité, universalité– dans l'oeuvre grothendieckienne, (2) les références directes de Grothendieck à Riemann dans la période 1949-1957 et les échos riemanniens dans *Récoltes et semailles* (1986), (3) les exemples de variable complexe de la première décennie 1949-1957 et les idées associées, (4) la vision émergeante au Congrès International de Mathématiciens en 1958, (5) quelques développements postérieures dans la deuxième (1958-1970) et la troisième (1981-1991) décennies des travaux connus de Grothendieck[4].

**1. Souplesse, naturalité, universalité**

On peut caractériser Grothendieck comme *le* mathématicien par excellence à la recherche permanente de *(i)* souplesse, *(ii)* naturalité, *(iii)* universalité. La souplesse *(i)* se réalise autour des formes simples, des structures lisses, des faisceaux; la naturalité *(ii)* est l'apanage des définitions catégoriques, des transformations naturelles, de la "mer montante"; l'universalité *(iii)* s'obtient autour de toutes sortes de théorèmes d'unification et de projectivité, avec l'émergence d'archétypes globaux qui gouvernent divers types locaux. Dans ce sens, les paroles de Grothendieck sont éclairantes : "Ma principale motivation a été et reste celle de développer les outils algébriques d'une *généralité* et d'une *souplesse* suffisantes pour le développement de cette géométrie arithmétique encore dans sa prime enfance (...) *coeur vraiment de mes amours* avec la mathématique"[5].

Dans son intervention au *Séminaire de Lectures Grothendieckiennes* (9 janvier 2018), Laurent Lafforgue parlait d'une "qualité de vérité" au coeur des démarches de Grothendieck, et comment les "vérités" se posaient au-dessus des "axiomes". En un certain sens, ceci correspond à imaginer une série d'*archétypes* (vérités) qui gouvernent les *types* (axiomes), et qui, grâce à leur même *abstraction*, permettent de construire des transits souples, des corrélations naturelles, des treillis universels. En fait, c'est justement en *échappant* du spécifique et du concret, que, dans la généralité, diverses obstructions peuvent être éliminées. La force de l'abstraction consiste précisément en permettre une *haute* architecture –souple, naturelle, universelle– qui descend après vers le *bas*, avec la résolution de problèmes particuliers[6]. Le *back-and-forth* entre abstraction et concrétion, entre universalité et particularité, méthode grothendieckienne s'il en est, capture la souplesse de la pensée mathématique.

---

[4] Par contre, tout est encore à connaître des 50.000 pages laissées à sa mort. Tout au moins, dans les quelques milliers de pages dédiées à la physique (selon rapport de Maltsiniotis), il est probable que la variable complexe fasse son apparition.
[5] Grothendieck, Alexandre (1985-1986), *Récoltes et semailles*, manuscrit, partie 4, page 1207 (nos accents).
[6] L'utilisation du topos étale d'un schéma pour résoudre les conjectures de Weil, ou l'utilisation du topos arithmétique (Connes) pour s'approcher à une résolution de l'hypothèse de Riemann, sont des frappants exemples d'une magnifique utilisation des machineries abstraites –souples, naturelles, universelles– avant de revenir sur les situations concrètes.

En fait, c'est grâce à une imagination sans bornes, dans l'espace des possibilités pures, que la mathématique développe ses outils et ses exemples les plus souples. Les surfaces de Riemann[7], le théorème de Riemann-Roch[8], les faisceaux algébriques cohérents[9], sont des telles constructions, fondamentales pour Grothendieck dans sa première décennie (1949-1957).

**2. Références directes à Riemann (1949-57 et *Récoltes et semailles*)**

Dès ses débuts à Bourbaki, Grothendieck rencontre Riemann. En effet, étant devenu cobaye en 1950, l'année suivante il fait une exposition de Riemann-Roch selon Kodaira (*Tribu 1951*). Comme nous avons signalé, l'influence de Riemann-Roch sera centrale dans son oeuvre. Entre 1953 et 1956, les références directes à la *sphère de Riemann* sont importantes : *(i) espace nucléaire* des fonctions holomorphes sur la sphère[10], *(ii) dualité* dans les espaces de fonctions holomorphes sur la sphère[11], *(iii) groupe structurel* d'un fibré holomorphe sur la sphère (somme de fibrés de fibre le plan complexe)[12], *(iv) catégorie additive, non abélienne*, donnée par les espaces de fibrés holomorphes sur la sphère[13].

---

[7] Les *surfaces de Riemann* apparaissent dans sa Thèse (1851), ibid., §5. Étant donnée une fonction de la variable complexe, la fonction n'est souvent pas injective (multiple-à-un), et il se pose le problème de son *inverse*. Par exemple, les fonctions puissance $z^n$ sont des fonctions ($n$-à-1) (un cadrant angulaire de taille $2\pi/n$ recouvre tout le plan complexe), la fonction exponentielle $e^z$ est ($\infty$-à-1) (une bande de hauteur $2\pi$ recouvre tout le plan). La façon usuelle d'inverser est de *restreindre localement le domaine* de la fonction. Riemann a imaginé la méthode opposée : *élargir globalement le codomaine*. On obtient des nouveaux objets géométriques globaux, construits en *collant* diverses feuilles locales : ce sont les surfaces de Riemann, sur lesquelles les inversions globales des fonctions sont possibles.

[8] Riemann-Roch (1857) unifie deux approches très différentes, en vue d'obtenir un *invariant naturel* pour une surface complexe. D'un côté, le *genre g* de la surface est obtenu en comptant le nombre de trous de la surface, ou, de façon équivalente, le nombre de coupes (moins un) avec lequel la surface devient disconnexe (*e.g.* le genre de la sphère est 0, le genre du tore 1, etc.) D'un autre côté, on peut penser aux "bonnes" fonctions (holomorphes) et aux "mauvaises" fonctions (méromorphes) qui peuvent être définies sur la surface. Si on fixe $n$ points sur la surface, on peut considérer l'espace (vectoriel) *Hol* des fonctions holomorphes avec des zéros sur ces points, et l'espace (vectoriel) *Mer* des fonctions méromorphes avec des pôles sur ces points. Le *théorème de Riemann-Roch* dit que $n - g + 1 = \dim(Mer) - \dim(Hol)$. De cette façon, un invariant géométrique *intrinsèque* (le genre) est relié à un invariant différentiel *extrinsèque* (la "différence harmonique" entre *Mer* et *Hol*). Au coeur de connections profondes entre géométrie, topologie, variable complexe, géométrie différentielle, algèbre et théorie des nombres, Riemann-Roch marque les débuts de la géométrie algébrique moderne. Son impact sera très fort chez Grothendieck, justement grâce à son caractère *rayonnant* sur des larges régions des mathématiques.

[9] Étant donné un faisceau dont ses fibres sont des anneaux, on peut itérer la construction et considérer le faisceau sur l'espace de ces anneaux, avec des fibres données par des modules sur ces anneaux. Si les transitions structurelles verticales fonctionnent bien (modules de type fini) et si les transitions homomorphes horizontales le font aussi (noyaux de type fini), nous avons un *faisceau algébrique cohérent* (Serre, Jean-Pierre (1955), "Faisceaux algébriques cohérents", in : *The Annals of Mathematics* 61.2, pp. 197-278). La *souplesse* est complète, donnée par les meilleurs conditions possibles (type fini) de la situation. Le faisceau des germes de fonctions holomorphes est un faisceau cohérent (Oka 1950).

[10] Grothendieck, Alexandre (1949-1953), *Produits tensoriels topologiques et espaces nucléaires*, Providence : American Mathematical Society, partie 2, page 58.

[11] Grothendieck, Alexandre (1953a), "Sur certains espaces de fonctions holomorphes I", *Journal für die reine und angewandte Mathematik* 192 : 35-64, §4. Le lieu central de la dualité dans la pensée grothendieckienne a été soulignée dans les divers exposés au *Séminaire de Lectures Grothendieckiennes* : Cartier (24 octobre 2017), Connes (7 novembre 2017), Lafforgue (9 janvier 2018), Pisier (6 février 2018), Caramello (6 mars 2018).

[12] Grothendieck, Alexandre (1955), "Sur la classification des fibrés holomorphes sur la sphère de Riemann", *American Journal of Mathematics* 79 : 121-138, p. 122.

[13] Grothendieck, Alexandre (1955-56), "Sur quelques points d'algèbre homologique", *Tôhoku Math. Journal* 9 : 119-221, p. 127.

D'un autre côté, la référence centrale à Riemann dans cette période se situe autour du théorème de Riemann-Roch généralisé[14]. C'est un excellent exemple des avantages que porte l'*orientation abstraite* de Grothendieck. Pour une variété *X* avec suffisamment de conditions de souplesse, Grothendieck imagine le *groupe libre K(X)* engendré para *tous* les faisceaux cohérents sur *X*. Cohérence, liberté, totalité –formes diverses de simplicité– aident à comprendre pourquoi *K(X) pourrait fonctionner* comme un archétype. Le théorème central de Grothendieck dans sa *K*-théorie démontre que, en effet, il *fonctionne* comme tel : *K* devient un *foncteur* relié au foncteur homologique *H*, à travers une transformation naturelle *C* donnée par les classes de Chern[15]. À ce moment-là, une obstruction à la commutativité émerge (*CK* n'est pas égal à *HC*), mais elle peut-être factorisée grâce à des classes de Todd[16], ce qui produit une formule généralisée Riemann-Roch-Serre-Hirzebruch. Le cas particulier de Riemann-Roch s'obtient en spécifiant sur une variété réduite à un point. Ainsi, l'archétype général (construction catégorique universelle *K(X)*), lorsqu'il est projeté sur une homologie triviale, capture le type spécifique (l'équation concrète de la surface, $n - g + 1 = \dim(Mer) - \dim(Hol)$).

Les références à Riemann deviendront de plus en plus nombreuses dans les travaux postérieurs de Grothendieck (pour un bref aperçu, voir notre *Section 5*, en bas). En particulier, Riemann apparaît directement dans plusieurs passages de *Récoltes et semailles* : *(i)* comme un des Maîtres dont il aspire à se ranger dans leur lignée[17], *(ii)* comme moteur d'un des thèmes centraux de son oeuvre[18], *(iii)* comme penseur révolutionnaire de la dualité continu/discret[19], *(iv)* comme rénovateur des mathématiques[20], *(v)* comme découvreur de profondes structures algébriques derrière la variable complexe[21]. En résumé, nous voyons Riemann cité à des moments cruciaux, quand Grothendieck cherche à *sonder les profondeurs*[22] de la pensée mathématique.

---

[14] Grothendieck, Alexandre (1957a), "Classes de faisceaux et théorème de Riemann-Roch" (Rapport à Serre 1957), reproduit dans *SGA*, vol. 6, pp. 20-77.

[15] Grothendieck, Alexandre (1957a), pp. 40-41, 63-64.

[16] Grothendieck, Alexandre (1957a), p. 72.

[17] "Moi qui ne suis pas fort en histoire, si je devais donner des noms de mathématiciens dans cette lignée-là, il me vient spontanément ceux de Galois et de Riemann (au siècle dernier) et celui de Hilbert (au début du présent siècle)", *Récoltes et semailles*, Préface, p. 13.

[18] "Yoga Riemann-Roch-Grothendieck (*K*-théorie, ...)", *Récoltes et semailles*, Préface, p. 21.

[19] "Il doit y avoir déjà quinze ou vingt ans, en feuilletant le modeste volume constituant l'oeuvre complète de Riemann, j'avais été frappé para une remarque de lui en passant. Il y fait observer qu'il se pourrait bien que la structure ultime de l'espace soit *discrète*, et que les représentations *continues* que nous nous en faisons constituent peut-être une simplification (excessive peut-être, à la longue...) d'une réalité plus complexe", *Récoltes et semailles*, Préface, p. 58.

[20] "L'oeuvre de Riemann (1826-1866) tient en un modeste volume d'une dizaine de travaux (il est vrai qu'il est mort dans la quarantaine), dont la plupart contiennent des idées simples et essentielles qui ont profondément renouvelé la mathématique de son temps", *Récoltes et semailles*, Partie 2, p. 227.

[21] "(...) les passages au quotient avaient un sens précis, dans le domaine analytique complexe, et les théorèmes à la Riemann-Serre (type GAGA) assuraient que le quotient final (qui était une courbe complexe compacte) avait une structure canonique de courbe *algébrique*", *Récoltes et semailles*, Partie 4, p. 1132.

[22] Grothendieck est ici proche de *Moby-Dick* (1851), grand roman des profondeurs de l'âme. Selon des témoignages de Tate (John et Karin) et de Yamashita, nous savons que Grothendieck aimait beaucoup le roman, et en lisait souvent des fragments.

## 3. Exemples et idées autour de la variable complexe (1949-1957)

Dans la table suivante, nous présentons des références, sinon exhaustives, au moins très représentatives, des occurrences de la variable complexe (formes d'holomorphie) dans la première décennie (1949-1957) des travaux de Grothendieck.

| *oeuvre*[23] | *pages, paragraphes* | *exemples* |
|---|---|---|
| 1949-1953 | 2e partie, p. 38 | (1) espace nucléaire de fonctions holomorphes |
| 1953a | § 4 | (2) dualité dans les espaces de fonctions holomorphes |
| 1953b | *passim* | (3) topologie (EVT) des espaces de fonctions holomorphes |
| 1953c | pp. 18, 73 | (4) normes complexes, constante de Grothendieck |
| 1955 | p. 122 | (5) groupe structurel d'un fibré holomorphe |
| 1955-1956 | pp. 127, 144 | (6) fibrés holomorphes, germes de fonctions holomorphes |
| 1957a | *passim* | (7) *K*-théorie : généralisation de Riemann-Roch |
| 1957b | §§ 6-8 | (8) faisceaux analytiques cohérents |

Autour de ces divers exemples (1-8), nous voyons surgir quelques-unes des idées essentielles de l'oeuvre grothendieckienne, mentionnées plusieurs fois au *Séminaire de Lectures Grothendieckiennes* : (1, 3) *nucléarité* (Cartier, 24 octobre 2017); (2) *dualité* (Lafforgue, 9 janvier 2018); (4) *projectivité* (Pisier, 6 février 2018); (5, 6, 7) *schématisation* (Szczeciniarz, 5 décembre 2017); (6, 7, 8) *localisation* (Connes, 7 novembre 2017; Caramello, 6 mars 2018).

---

[23] En plus des références déjà données, nous mentionnons ici d'autres travaux : Grothendieck, Alexandre (1953b), "Sur certains espaces de fonctions holomorphes II", *Journal für die reine und angewandte Mathematik* 192 : 77-95; Grothendieck, Alexandre (1953c), "Résumé de la théorie métrique des produits tensoriels topologiques", *Bol. Soc. Mat. Sao Paulo* 8 (publié 1956) : 1-79; Grothendieck, Alexandre (1957b), "Sur les faisceaux algébriques et les faisceaux analytiques cohérents", *Séminaire Cartan ENS* 9 (2) : 1-16.

Le réseau d'idées qui émerge de ces travaux représente ainsi, en forme adéquate, la *complexité mathématique* de l'oeuvre de Grothendieck, soulignée maintes fois par tous les expositeurs au *Séminaire de Lectures Grothendieckiennes*. On peut résumer cette complexité dans un *va-et-vient permanent* entre types et archétypes, entre structures concrètes et formes générales, entre obstructions particulières et transits universaux. Dans ce va-et-vient, la variable complexe agit du côté des *types*, avec ses exemples plastiques et ses harmonies naturelles. Quand les structures de la variable complexe (espaces, fibrés, faisceaux) rentrent en contact avec des *archétypes* universaux (constante de Grothendieck, catégories abéliennes, *K*-théorie), on se plonge dans le *coeur profond* des mathématiques ("coeur vraiment de mes amours avec la mathématique"). Les catégories abstraites, générales, universelles, agissent comme des "mers montantes" qui dissolvent les difficultés. Au niveau des *archétypes-vérités* (Lafforgue, 9 janvier 2018), une grande *souplesse* gouverne la créativité mathématique. Il est frappant d'observer que cette souplesse reproduit en grande mesure, à un niveau abstrait supérieur, l'élasticité inhérente des formes holomorphes et méromorphes, au niveau concret inférieur de la variable complexe.

## 4. ICM 1958

Promu en jeune mathématicien vedette avec sa preuve de Riemann-Roch, Grothendieck est invité à donner une conférence plénière au Congrès International de Mathématiciens (ICM) en 1958. Avec trente ans à peine, il a déjà renouvelé de fond en comble l'analyse fonctionnelle[24], développé l'analyse fine des espaces de Banach[25], donné les bases de la théorie infinitaire des catégories et ses connexions avec l'homologie[26], et inventé la *K*-théorie[27], offrant des labeurs pour des générations entières de mathématiciens. Et soudain, au ICM 1958, sans un moindre indice au préalable, il imagine les schémas et les topos, introduit les premiers, et donne en quelques pages le dessin de son gigantesque programme[28] pour démontrer les conjectures de Weil. Le *coeur* de son idée centrale se réduit à cinq lignes (!):

> Although interesting relations must certainly exist between these cohomology groups and the "true ones", it seems certain now that the Weil cohomology has to be defined by a completely different approach. Such an approach was recently suggested to me by the

---

[24] Grothendieck (1949-1953), *Produits tensoriels topologiques et espaces nucléaires*, op. cit.
[25] Grothendieck (1953c), "Résumé de la théorie métrique des produits tensoriels topologiques", op. cit.
[26] Grothendieck (1955-56), "Sur quelques points d'algèbre homologique", op. cit.
[27] Grothendieck (1957a), "Classes de faisceaux et théorème de Riemann-Roch", op. cit.
[28] Le programme a été développé à l'IHES, entre 1958 et 1970, dans ses *Éléments de Géométrie Algébrique* (avec Dieudonné) (1959-1964) et son *Séminaire de Géométrie Algébrique* (avec une pléiade de brillants élèves) (1960-1969).

*connections between sheaf-theoretic cohomology and cohomology of Galois groups on the
one hand, and the classification of unramified coverings of a variety on the other* (...)[29]

De cette façon, à la recherche d'une *vraie* cohomologie[30], Galois, d'un côté, et Riemann, de l'autre, sont connectés dans ce que deviendra la notion de *schéma étale*[31], clef technique fondamentale du programme. On voit bien comment la *variable complexe ronde* derrière la vision géniale "connections between sheaf-theoretic cohomology and cohomology of Galois groups on the one hand, and the classification of unramified coverings of a variety on the other". D'un côté, les faisceaux, provenant de la topologie (Leray) et de la variable complexe (Cartan), sont imbriqués avec les travaux de Galois sur les extensions des corps de nombres dans la clôture algébrique des nombres complexes; d'un autre côté, les recouvrements non ramifiés d'une variété, provenant des voisinages lisses et de la continuation analytique dans la variable complexe, sont imbriqués avec les recouvrements universels d'une surface de Riemann. C'est en *élevant* Galois, Riemann et la variable complexe au monde abstrait des schémas et des topos, que Grothendieck découvre la *souplesse* suffisante pour arriver à la "vraie" cohomologie nécessaire pour résoudre les conjectures de Weil.

## 5. Quelques développements postérieurs (1958-1970, 1981-1991)

Autour de 1955, une table (voir *Figure 1*) permet d'observer les diverses forces en action dans la pensée grothendieckienne. Comme on peut le voir, suivant notre "calibrage des forces" dans la figure, la place de la variable complexe est, sinon centrale, du moins bien importante. D'autre part, les thèmes de la variable complexe, et les développements autour de Riemann en particulier, apparaissent nombre de fois dans les écrits postérieurs (1958-1970 & 1981-1991) de Grothendieck. Il serait impossible d'en décrire ici les divers chemins, mais une courte liste de telles apparitions est la suivante: *(i)* travaux de fondements de la géométrie analytique (comprise, à la Serre, comme géométrie des fonctions analytiques)[32], *(ii)* étude de diverses propriétés de souplesse dans les constructions de la géométrie algébrique[33], *(iii)* exploration des espaces moduli, groupe de Grothendieck-Teichmüller et géométrie anabélienne[34], *(iv)* vision nouvelle des surfaces de Riemann autour des dessins d'enfants[35], *(v)* reconstruction de l'homotopie *at large*, en étudiant les propriétés fonctorielles de constructions catégoriques abstraites qui modèlent les CW-complexes (homologie, groupes complexes)[36]. On rencontre ainsi une *profonde continuité* dans les divers travaux de Grothendieck autour de la variable complexe, ce qui, maintenant, ne devrait plus nous surprendre: la *souplesse*, la *naturalité*, l'*universalité* de la variable complexe se reflètent parfaitement sur l'ambiance abstraite des catégories, où, suivant le vers de Baudelaire, "tout est luxe, calme et volupté".

---

[29] Grothendieck, Alexandre (1958), "The cohomology theory of abstract algebraic varieties", dans: *Proceedings International Congress of Mathematicians (Edinburgh)*, Cambridge: Cambridge University Press, pp. 103-118, citation p. 104 (phrases accentuées par Grothendieck lui-même).

[30] Les échos de la conférence de Lafforgue (9 janvier 2018) au *Séminaire de Lectures Grothendieckiennes* sont particulièrement pénétrants ici. Grothendieck cherche une *vérité archétypique*, bien au-delà des types donnés par des structures et des axiomes.

[31] Un *schéma affine* est un faisceau de la forme "faisceau de représentation d'un anneau" (étant donné un anneau *A* commutatif unitaire, son spectre *S* est formé par les idéaux premiers de *A*, et on peut mettre sur *S* la topologie de Zariski; après, sur chaque idéal *P* du spectre, on prend comme fibre sa localisation $A(A-P)^{-1}$; en unissant toutes les fibres, et en mettant la topologie initiale sur l'union, on obtient un faisceau sur *S*: c'est le faisceau de représentation de *A*). Un *schéma* est un faisceau qui localement est isomorphe à un schéma affine (en d'autres termes, un schéma est un *recollement* de schémas affines). Un morphisme entre schémas est *étale* s'il se comporte de la façon la plus *souple* possible "en regard de Galois et Riemann", c'est-à-dire, s'il est de type fini et plat (Galois) et non ramifié (Riemann) (*SGA*, 1, p. 4). Étant donné un schéma *X*, la collection de *tous* les morphismes étales sur *X* forme une topologie de Grothendieck (*SGA* 4.1, p. 219), à partir de laquelle on construit le *topos étale d'un schéma* (*SGA*, 4.2, p. 343), qui à son tour donne lieu naturellement à la *cohomologie étale*, clef de la preuve des conjectures de Weil.

[32] Grothendieck, Alexandre (1961), "Techniques de construction en géométrie analytique", *Séminaire Cartan* 13, exposés 7-17, Paris: Secrétariat Mathématique ENS.

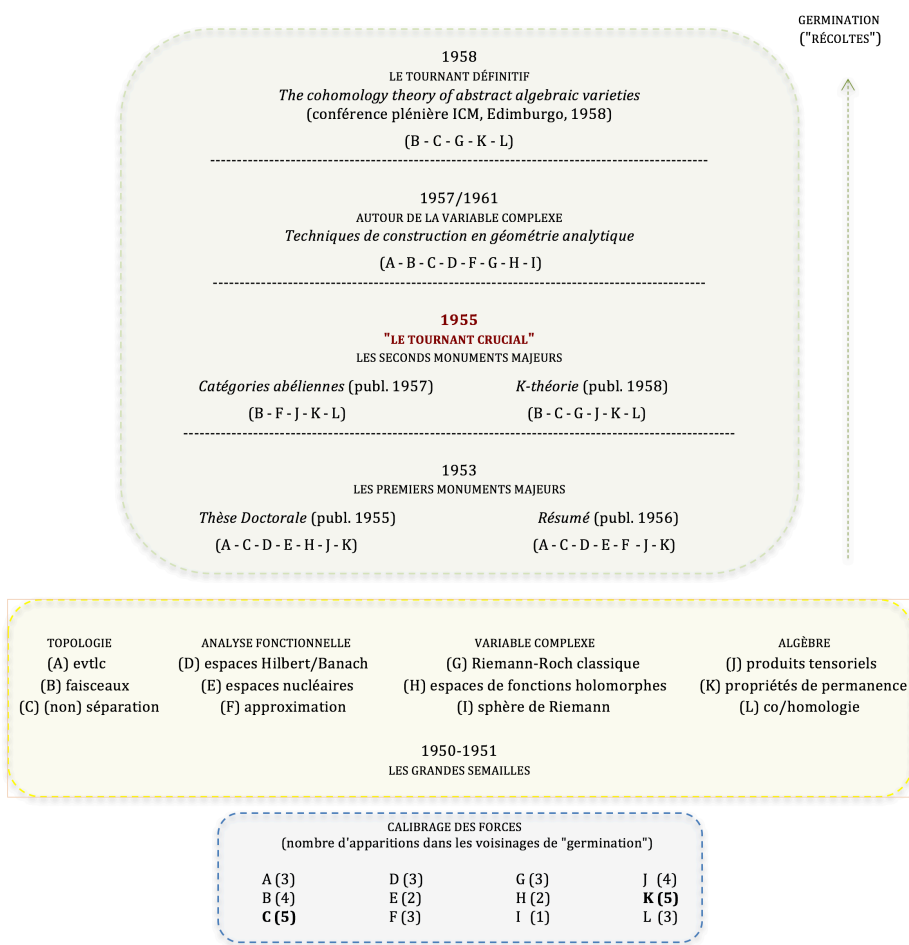

*Figure 1. La pensée de Grothendieck autour de 1955*

---

[33] Grothendieck, Alexandre & Dieudonné, Jean (1959-1964), *Éléments de Géométrie Algébrique*, Paris: IHES; Grothendieck, Alexandre, *et. al.* (1960-1969), *Séminaire de Géométrie Algébrique*, Paris: IHES.
[34] Grothendieck, Alexandre (1981), "La longue marche à travers la théorie de Galois", manuscrit.
[35] Grothendieck, Alexandre (1984), "Esquisse d'un programme", manuscrit.
[36] Grothendieck, Alexandre (1991), "Les dérivateurs", manuscrit.

**Conclusion.**

La pensée mathématique de Grothendieck agit à plusieurs niveaux, suivant un incessant *back-and-forth* entre le concret (exemples particuliers) et l'abstrait (formes universelles). La *réduction* de ce va-et-vient à des niveaux purement généraux ("abstract nonsense") a montré, chez divers commentateurs, une lecture très pauvre (ou plutôt, une non lecture) des écrits de Grothendieck. Bien au contraire, ceux-ci se trouvent *remplis* de toutes sortes d'exemples et de considérations précises sur des régions mathématiques traditionnelles (analyse fonctionnelle, topologie, algèbre abstraite, théorie des nombres, combinatoire, géométrie élémentaire, etc.) Dans cette vision du concret, la variable complexe joue un rôle déterminant. La *continuité* du regard grothendieckien sur la variable complexe –de la *Thèse* (1949-53) jusqu'à l'*Esquisse d'un programme* (1984) et *Les dérivateurs* (1991)– montre l'importance de profiter de cette région *locale* des mathématiques pour raffiner une vision *globale*. En fait, on peut facilement imaginer une sorte de *faisceau métamathématique* dans son oeuvre, partant du local complexe vers le global catégorique. Dans ce long parcours, Grothendieck restera toujours fidèle au coeur vraiment de ses amours, une géométrie arithmétique imprégnée de perspectives de la variable complexe et de la topologie.